\documentclass[aps,preprint,superscriptaddress,floatfix,amsmath,amssymb]{revtex4}
\usepackage{inputenc}
\usepackage[a-1b]{pdfx}
\usepackage{amsmath}
\usepackage{amssymb}
\usepackage{amsthm}
\usepackage{algcompatible}
\usepackage{algpseudocode}
\usepackage{newfloat}
\usepackage{bm}
\usepackage{graphicx}
\usepackage{subfig}
\usepackage{hyperref}
\usepackage[capitalise]{cleveref}

\DeclareFloatingEnvironment[
    fileext=loa,
    listname=List of Algorithms,
    name=ALGORITHM,
    placement=htbp,
]{algorithm}

\begin{document}

\title{Approximation of Optimal Control Surfaces for $2\times 2$ Skew-Symmetric Evolutionary Game Dynamics}

\author{Gabriel Nicolosi}
\email{gxr286@psu.edu}
\author{Terry Friesz}
\email{tlf13@psu.edu}
\affiliation{Harold and Inge Marcus Dept.\ of Industrial and Manufacturing Engineering, Penn State University,
University Park, PA 16802}

\author{Christopher Griffin}%
\email{griffinch@psu.edu}
\affiliation{Applied Research Laboratory, Penn State University,
University Park, PA 16802}

\date{July 14, 2022}

\begin{abstract}
In this paper we study the problem of approximating the general solution to an optimal control problem whose dynamics arise from a $2\times 2$ skew-symmetric evolutionary game with arbitrary initial condition. Our approach uses a Fourier approximation method and generalizes prior work in the use of orthogonal function approximation for optimal control. At the same time we cast the fitting problem in the context of a non-standard feedforward neural network and derive the back-propagation operator in this context. An example of the efficacy of this approach is provided and generalizations are discussed.
\end{abstract}

\maketitle


\section{Introduction}
\label{sec:Introduction}

Optimal control and variational problems have been studied extensively, see \cite{K04,F10,SL12} among numerous other texts. Control of evolutionary game dynamics is a relatively recent problem in the area of non-linear dynamics and non-linear control. In \cite{PQ12} Pantoja and Quijano investigate a distributed optimization problem on a network with the replicator. More recently \cite{GTBS16} studies convergence of best-response strategies on graphs. Fan and Griffin \cite{GF22} study optimal control of odd circulant games generalizing work in \cite{FG17} in which control of the Bass model is considered. At the same time, there has been extensive work on reinforcement learning (RL) based methods for control of dynamical systems \cite{KLM96,KBP13} with more recent work in (deep) neural network based methods coming to the fore \cite{ADBB17}. Optimal control is a natural component of this broader area of RL research \cite{KVML17}. 

In this paper, we consider the problem of optimal control of evolutionary game dynamics through the lens of Fourier analysis. We note this area has been widely studied in control theory \cite{CS87,RTA89,NY90,RR90,YN91,T91,YC94,MK05} and has found applications in physics, chemistry and materials science \cite{HMM00,KHK10}. The general area falls under the use of orthogonal functions in \textcolor{black}{the context of direct collocation methods for trajectory optimization and optimal control \cite{R79,R14,K17}}. \textcolor{black}{Orthogonal functions in systems and control are} summarized in \cite{DM95} with more recent work by Ragazzi focusing on Legendre polynomials rather than trigonometric polynomials  \cite{RHF95,RH96,RS98,RY02,MR03,MR04}. The previous methods focus on methods of integration using orthogonal polynomial methods for solving the Riccati equations and two-point boundary value problems that emerge as a result of optimal control problems. 

We vary this approach inspired both by the Fourier methods and methods for approximating solutions to the corresponding Hamilton-Jacobi-Bellman (HJB) equations that arise from optimal control problems \cite{M07,ACC08,KK18,JP20,GSS21}. These methods attempt to approximate the value function of the HJB equation and then use it to construct an optimal controller. In contrast, in this paper we consider an optimal control problem of Lagrange type
\begin{equation}
\begin{aligned}
\min_{u} \;\; & \int_0^T f(x,u) \, dt\\
s.t.\;\; & \dot{x} = g(x,u)\\
& x(0) = x_0\\
& u \in L_2([0,T]).
\end{aligned}
\label{eqn:OptCon}
\end{equation}
where $f:\mathbb{R} \times L_2([0,T]) \rightarrow \mathbb{R}$, $x \in X \subseteq \mathbb{R}$. In what follows, we will assume that $g(x,u)$ \textcolor{black}{will be constructed} from a two-strategy skew-symmetric evolutionary game. For fixed $x_0$, the \textit{open loop} optimal controller is the function $u(t)$ solving \cref{eqn:OptCon}. In this paper, we focus specifically on the open-loop optimal control problem, leaving the closed-loop control problem, \textcolor{black}{i.e., a control of the type $u=u(x,t,x_0)$}, for future work. If $\textcolor{black}{x(0) =} x_0$ is unknown a priori \textcolor{black}{and $t \in [0,L]$}, then our objective is to approximate an optimal control surface $u(t,x_0)$ so that for fixed $x_0$, $u(t,x_0)$ is the optimal open-loop controller given the fixed initial condition $x_0$. That is, instead of approximating the solution to the non-linear HJB equation and using this approximation to construct an optimal control for arbitrary initial condition, we approximate the optimal control surface directly. The main contributions of this paper are:
\begin{enumerate}

\item We extend the work in \cite{CS87,RTA89,NY90,RR90,YN91,T91,YC94,MK05} to approximate not only the optimal control for a fixed initial condition but for an arbitrary initial condition making the approach more like the approximation to the HJB \cite{JP20}.

\item Unlike work in \cite{NY90}, which specifically eschews a gradient based method, we construct an explicit gradient descent method that can be used like back-propagation in a neural network. 

\item The approach to estimating an optimal control surface is applied to non-linear dynamics arising from a two-strategy evolutionary game where we show excellent performance for reasonably small size approximations.

\end{enumerate} 
The remainder of this paper is organized as follows: In \cref{sec:ProblemStatement} we layout the proposed structure of the optimal control surface $u(t,x_0)$ and show its relation to a non-standard neural network problem. We also discuss the assumed evolutionary game dynamics that govern the state equation. Construction of the back-propagation operator is provided in \cref{sec:BackPropagation} along with the optimization algorithm for approximating the optimal control surface. In \cref{sec:ExperimentalResults} we provide experimental results. Generalizations are discussed in \cref{sec:Generalization}. Conclusions and future directions are provided in \cref{sec:Conclusion}.

\section{Problem Construction and Preliminaries}\label{sec:ProblemStatement}
\begin{figure}[htbp]
\centering
\includegraphics[width=0.65\textwidth]{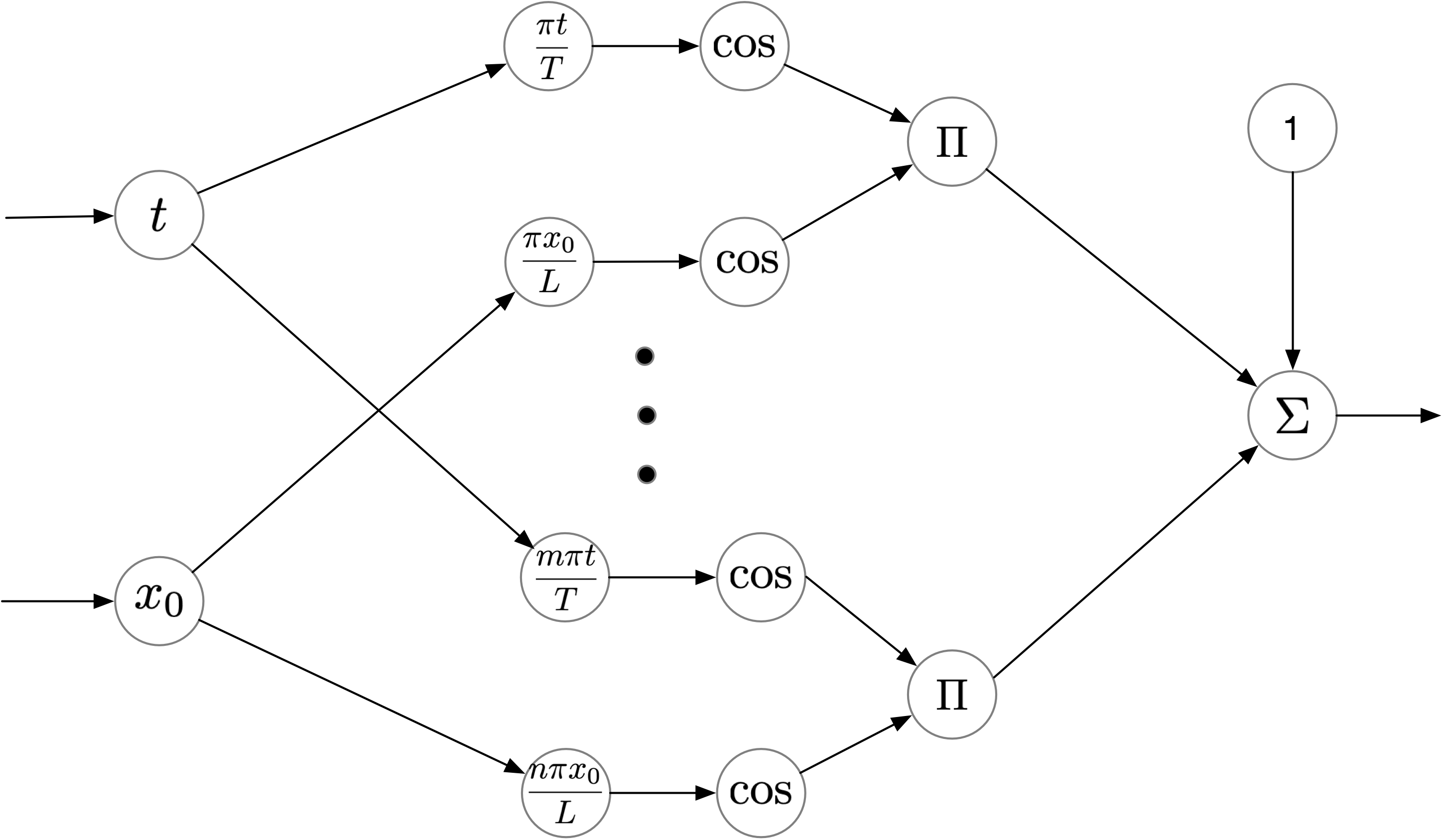}
\caption{The formulation of the optimal control surface as a non-standard feedforward neural network inspires the use of the construction of the back-propagation operator.}
\label{fig:NN}
\end{figure}

We assume an approximation of $u(t,x_0)$ as
\begin{equation}
u = u(t,x_0) \approx \sum_{m=0}^M \sum_{n=0}^N a_{mn} \cos\left(\frac{m\pi t}{T}\right) \cos\left(\frac{n\pi x_0}{L}\right).
\label{eqn:ufourier}
\end{equation}
Using this construction, our problem reduces to identifying the finite set of Fourier coefficients of true (hidden)  optimal control surface $u(t,x_0)$. Unlike an ordinary Fourier approximation, our goal is not to find $u(t,x_0)$ and then build $a_{mn}$ but rather to build the coefficients directly from the primal problem \cref{eqn:OptCon}. This structure can be represented as a non-standard feedforward neural network (\cref{fig:NN}). We note that optimal control problems have been addressed using standard neural network architectures in \cite{EP13} with some success. 
Phrasing this approximation problem in the context of a non-standard neural network inspires our construction of a back-propagation operator.

Suppose we have a finite sample $X_0$ of initial conditions. Let $\mathbf{a}$ be a vector composed of the Fourier coefficients. Then define \textcolor{black}{
\begin{align}
    J(x_0,\mathbf{a}) =  &\int_0^T  f(x,u) \, dt\\
    s.t.\;\; &\dot{x} = g(x,u) \label{eqn: stdynamics}\\
    &x(0) = x_0 \label{eqn: stinitialcond}\\
    &u(t,x_0) \approx \sum_{m=0}^M \sum_{n=0}^N a_{mn} \cos\left(\frac{m\pi t}{T}\right) \cos\left(\frac{n\pi x_0}{L}\right). \label{eqn: stcontrolapprox}
\end{align}
The objective function and its constraints are then given by
\begin{equation}
J(\mathbf{a}) = \sum_{x_0 \in X_0} J(x_0;\mathbf{a}) \quad s.t. \quad\text{\cref{eqn: stdynamics,eqn: stinitialcond,eqn: stcontrolapprox} hold} \; \forall \; x_0 \in X_0.
\label{eqn: objfunc}
\end{equation}
We have now converted the time continuous problem in \cref{eqn:OptCon} into a nonlinear programming problem with finite decision variables
\begin{equation}
\min_\mathbf{a} \;\; J(\mathbf{a}). 
\label{eqn:DiscreteProblem}
\end{equation}}
Problem \cref{eqn:DiscreteProblem} can be solved using a direct optimization technique (e.g., LGBFS, \textcolor{black}{conjugate gradient,} etc.) However, the remainder of this paper will be dedicated to casting this into a non-standard neural network architecture for a specific class of optimal control problems and then constructing the back-propagation operator for this neural network structure.

\subsection{Evolutionary Game Dynamics}
In deriving a back-propagation operation to solve \cref{eqn:DiscreteProblem}, we assume the equations of motion are given by the replicator dynamics. Let $\mathbf{A} \in \mathbb{R}^{n \times n}$ be a (payoff) matrix. The replicator dynamics are given by
\begin{displaymath}
\dot{x}_i = x_i\left(\mathbf{e}_i - \mathbf{x}\right)^T\mathbf{A}\mathbf{x},
\end{displaymath}
where $\mathbf{x} = \langle{x_1,\dots,x_n}\rangle \in \Delta_{n-1}$ \textcolor{black}{and $\mathbf{e}_i$ is the $i^{th}$ standard unit vector} . Here $\Delta_{n-1}$ is the unit simplex.  We focus on the case when
\begin{displaymath}
\mathbf{A} = \begin{bmatrix} 0 & -\rho \\ \rho & 0\end{bmatrix}.
\end{displaymath}
The resulting replicator dynamics
\begin{align*}
    \dot{x}_1 = -\rho x_1 x_2\\
    \dot{x_2} = \rho x_1x_2
\end{align*}
describes rumor spreading \cite{ALM19}, \textcolor{black}{susceptible-infected (SI) epidemic dynamics} \cite{H00}, the aspatial component of Fisher's equation \cite{F37} and the Bass model of social science \cite{B69}. 
\textcolor{black}{Game theoretic analysis of epidemic dynamics have been studied previously from a game-theoretic and socio-physics context \cite{T15,T19,T21}. The donor/recipient approach studied by Tanimoto has particular relevance to this problem \cite{T09}. In this paper, we focus on the problem of controlling the trajectories of evolutionary game dynamics when the pay-off matrix itself is being manipulated. Recently this approach has found interest in the biomedical community \cite{ZCBG17,GSV20}, where empirical methods are being used to produce control strategies. Since this problem is inherently non-linear this paper develops an approximation method for the optimal controller that while simultaneously reducing the infinite dimensional optimization problem to a finite dimensional problem in the spirit of \cite{CS87,RTA89,NY90,RR90,YN91,T91,YC94,MK05}. We illustrate this approach on a two-strategy game because we can simplify the dynamics  by setting} $x_2 = 1 - x_1$. We can replace the two equations of the replicator dynamics with the single equation of motion
\begin{displaymath}
\dot{x} = \rho x(1-x).
\end{displaymath}
\textcolor{black}{More general classes of problems are discussed in future work.} If we assume parameter \textcolor{black}{$\rho$} in the payoff matrix is a (linear) function of the control $u$ so that
\begin{equation}
\rho = \beta u - \xi,
\end{equation}
then the equations of motion become
\begin{equation}
\dot{x} = x(1-x)\left(\beta u - \xi \right). \label{eqn: Dynamics}
\end{equation}
The specific variation of \cref{eqn:OptCon} is then
\begin{equation}
\begin{aligned}
\min_{u} \;\; & \int_0^T  \frac{k_1}{2}x^2 + Rx u + \frac{k_2}{2}u^2  \, dt\\
s.t.\;\; & \dot{x} = x(1-x)\left(\beta u - \xi \right)\\
& x(0) = x_0\\
& u \in L_2([0,T]).
\end{aligned}
\label{eqn:OptCon2}
\end{equation}
\textcolor{black}{With $R$, $k_1$ and $k_2$ being real valued constants. In \cref{eqn:OptCon2}, a quadratic objective function is suitably chosen to frame this problem in the general context of linear quadratic control problems. Even though this is the choice of objective function for this paper, the back-propagation operation defined in \cref{sec:BackPropagation} is extended to any other objective function as long as its derivative can be computed analytically.}

\section{Construction of the Back-Propagation Operation}
\label{sec:BackPropagation}
Back-propagation is simply a computational application of the chain rule combined with gradient descent \cite{GB16}. Our objective is to construct $\nabla J(\mathbf{a})$. Let $\varphi_{x_0}(t)$ be the flow satisfying the dynamics
\begin{equation}
\dot{x} = x(1-x)(\beta u - \xi)
\label{eqn:xdot}
\end{equation}
with $x(0) = x_0$. Then
\begin{equation}
J(\mathbf{a}) = \sum_{x_0 \in X_0} \int_0^T \frac{k_1}{2}\varphi_{x_0}(t)^2 + R\varphi_{x_0}(t)u + \frac{k_2}{2}u^2\,dt
\label{eqn: objectivefunctional}
\end{equation}
with $u$ given by \cref{eqn:ufourier}. Differentiating with respect to $a_{mn}$ yields
\begin{equation}
\frac{\partial J}{\partial a_{mn}} = \sum_{x_0 \in X_0} \left(\int_0^T
k_1\varphi_{x_0}\frac{\partial \varphi_{x_0}}{\partial a_{mn}} + R\left( u\frac{\partial \varphi_{x_0}}{\partial a_{mn}} + \varphi_{x_0}\frac{\partial u}{\partial a_{mn}}\right) + k_2u\frac{\partial u}{\partial a_{mn}}\,\, dt\right).
\end{equation}
Factoring yields
\begin{equation}
\frac{\partial J}{\partial a_{mn}} = \sum_{x_0 \in X_0} \left(\int_0^T \frac{\partial \varphi_{x_0}}{\partial a_{mn}}\left(k_1 \varphi_{x_0} + R u \right) + \frac{\partial u}{\partial a_{mn}}\left( R \varphi_{x_0} + k_2 u \right)
\,\,dt
\right).
\label{eqn:dJdamn}
\end{equation}
By assumption, $u$ is given by \cref{eqn:ufourier} and consequently
\begin{displaymath}
\frac{\partial u}{\partial a_{mn}} =\cos\left(\frac{m\pi t}{T}\right) \cos\left(\frac{n\pi x_0}{L}\right).
\end{displaymath}
An expression for $\varphi_{x_0}(t)$ can be obtained by integrating \cref{eqn: Dynamics}, where $u$ is given by the Fourier approximation in \cref{eqn:ufourier}
\begin{equation}
\int_0^t \frac{dx}{x(1-x)} = \int_0^t (\beta u(\tau) - \xi) \quad d\tau
\label{eqn: rhsintegration}
\end{equation}
We then obtain a closed-form expression for $\varphi_{x_0}$ in terms of the same coefficients $a_{mn}$ in $u$
\begin{equation}
    \varphi_{x_0} = \frac{1}{\left(1 + K_{x_0}\exp\left({V_{x_0}(t)}\right)\right)} 
\end{equation}
where $K_{x_0}$ given by 
\begin{equation}
K_{x_0} = \frac{1-x_0}{x_0}\exp\left({V_{x_0}(0)}\right)
\end{equation}
is the constant of integration in \cref{eqn: rhsintegration} and $V_{x_0}(t)$ is written as
\begin{equation}
    V_{x_0}(t) = \beta U_{x_0}(t) - \xi t
\end{equation}
in which $U_{x_0}(t)$ is the integral $\int_0^t u(\tau,x_0) d\tau$ of the control approximation given by \cref{eqn:ufourier}
\begin{equation}
    U_{x_0}(t) = \sum_{m=1}^M \sum_{n=0}^N a_{mn}\sin\left(\frac{m \pi t}{T}\right)\frac{T}{m \pi} \cos\left(\frac{n \pi x_0}{L}\right) + \sum_{n=0}^N a_{0n}\cos\left(\frac{n \pi x_0}{L}\right)t. 
\end{equation}
\textcolor{black}{By the chain rule}, an expression for $\frac{\partial \varphi_{x_0}}{\partial a_{mn}}$ \textcolor{black}{can be written as}
\begin{equation}
    \frac{\partial \varphi_{x_0}}{\partial a_{mn}} = \frac{\partial \varphi_{x_0}}{\partial V_{x_0}} \frac{\partial V_{x_0}}{\partial U_{x_0}} \frac{\partial U_{x_0}}{\partial a_{mn}},
    \label{eqn: dphidamn}
\end{equation}
where the right-hand side derivatives are computed as
\begin{align}
    \frac{\partial \varphi_{x_0}}{\partial a_{mn}} &= \frac{K_{x_0}\exp(V_{x_0}(t))}{(K_{x_0} + \exp(V_{x_0}(t)))^2} \label{eqn: chain1}\\
    \frac{\partial V_{x_0}}{\partial U_{x_0}} &= \beta \label{eqn: chain2} \\
    \frac{\partial U_{x_0}}{\partial a_{mn}} &= 
    \begin{cases}
    \sin\left(\frac{m \pi t}{T}\right) \cos\left(\frac{n \pi x_0}{L}\right) \frac{T}{m \pi} &\mbox{if } m \neq 0 \\
    \cos\left( \frac{n \pi x_0}{L} \right) t &\mbox{if } m = 0 \\
    \end{cases}
    \label{eqn: chain3}
\end{align}
Substituting \cref{eqn: dphidamn} into \cref{eqn:dJdamn} we obtain a closed-form expression for $\nabla J(\mathbf{a})$, allowing for the minimization of the functional $J(\mathbf{a})$ by performing gradient descent over the space spanned by the Fourier coefficients in $\mathbf{a}$.  

\subsection{Gradient Descent Over the Space of Fourier Coefficients}
A solution to \cref{eqn: objectivefunctional} is then made possible by performing the following procedure described in \cref{alg:cap}.
\begin{algorithm}
\caption{Minimization of $J(\mathbf{a})$ by Gradient Descent}\label{alg:cap}
\begin{algorithmic}
\Require Initialize the coefficients $\mathbf{a^0}$, initial control $u = u(\mathbf{a^0})$ and state $\textcolor{black}{\varphi} = \textcolor{black}{\varphi}(\mathbf{a^0})$, gradient tolerance $\epsilon$ and descent (learning) rate $\alpha$;
\State $k \gets 0$;
\State Compute $\nabla J(\mathbf{a^0})$;
\While{$\nabla J(\mathbf{a^k}) \ge \epsilon$}
    \State Update coefficients $\mathbf{a} \gets \mathbf{a} - \alpha \nabla J(\mathbf{a^k})$
    \State Update control approximation $u = u(\mathbf{a^k})$ and trajectory $\textcolor{black}{\varphi} = \textcolor{black}{\varphi}(\mathbf{a^k})$
    \State Compute $\nabla J(\mathbf{a^k})$
    \State $k \gets k + 1$
\EndWhile
\end{algorithmic}
\end{algorithm}
We note the parallels of \cref{alg:cap} with the back-propagation operation present in traditional neural networks \cite{RHW86}. In our case, instead of the traditionally employed loss functions (\textcolor{black}{e.g.} quadratic, logistic, \dots), the objective function being minimized is represented by the objective functional itself (\cref{eqn: objfunc}). This expression is implicitly parameterized by the Fourier coefficients in $\mathbf{a}$ through the approximated control and its integrated state trajectory arising from the dynamical (control dependent) state constraint. Thus, computing $\nabla J(\mathbf{a})$ requires the computation of $\frac{\partial u}{\partial a_{mn}}$ and $\frac{\partial \textcolor{black}{\varphi}}{\partial a_{mn}}$, the latter requiring a chain of intermediate derivatives (\cref{eqn: chain1,eqn: chain2,eqn: chain3}), similar to the traditional neural network setting wherein a chain of derivatives is constructed and the weights (coefficients) are adjusted through back-propagation.

\section{Experimental Results}
\label{sec:ExperimentalResults}
In this Section, we apply the procedure described above to obtain an approximated solution to \textcolor{black}{an example problem and discuss interpretable conclusions from the approximation. Consider a simplistic epidemic model in which the two strategies of the evolutionary game are susceptible and infected. In the absence of a treatment (intervention) the dynamics are given by:
\begin{align*}
\dot{x}_1 = -\xi x_1 x_2\\
\dot{x}_2 = \xi x_1 x_2,
\end{align*}
where $x_1$ is the proportion of the population that is susceptible and $x_2$ is the proportion of the population who are infected. The parameter $\xi$ is the standard infection rate. This is consistent with the work in \cite{T09,T15,T19,T21}. Intervention leads to the controlled dynamics:
\begin{align*}
\dot{x}_1 = (\beta u -\xi) x_1 x_2\\
\dot{x}_2 = (\xi-\beta u)x_1 x_2,
\end{align*}
where $u$ is a measure of the effort (input) made in providing the treatment. This yields the skew-symmetric evolutionary game dynamics given in \cref{eqn:OptCon2}. We assume the intervention has a quadratic cost $Cu^2$ (that should be minimized) and that societal benefit arises not only from $x_2$ being minimized but also from the interaction of the treatment effort $u$ with $x_1$. That is, high-impact treatment efforts will lead to improved public health outputs. Letting $x = x_1$ and $x_2 = 1 - x$, the optimal control problem to be considered is 
\begin{equation}
\begin{aligned}
\max_{u} \;\; & \int_0^T \alpha xu - Cu^2 \, dt\\
s.t.\;\; \dot{x} &= x(1-x)(\beta u-\xi)\\
x(0) &= x_0\\
u &\ge 0\\
\end{aligned}
\label{eqn:SoftwareMaintenanceProblem}
\end{equation} 
}

\textcolor{black}{In this problem $\alpha$ is the conversion factor measuring societal benefit from the interaction of individuals with the susceptible strategy and the treatment effort. It is interesting to note that the structure of this problem is similar to the software maintenance problem explored in \cite{FG17} and the cyber-bullying problem explored in \cite{QSGU18}. By exploiting the concavity  of \cref{eqn:SoftwareMaintenanceProblem}, an analytical solution can be obtained for the closed-loop controller \cite{FG17}. In contrast our approximation builds an open-loop controller for the evolutionary system.}


To solve \cref{eqn:SoftwareMaintenanceProblem} for a set of distinct initial conditions and find the (approximated) optimal control surface $u^*=u^*(t,x_0)$, we formulate \cref{eqn:SoftwareMaintenanceProblem} as in \cref{eqn:DiscreteProblem} and choose a set of distinct $x_0 \in X_0$ (approximation points). We assue the following parameter values: $\alpha = 2$, $C = 1$, $\beta=\frac{1}{2}$, $\xi = \frac{1}{4}$.

\begin{figure}[ht!]
\centering
\includegraphics[width=0.46\textwidth]{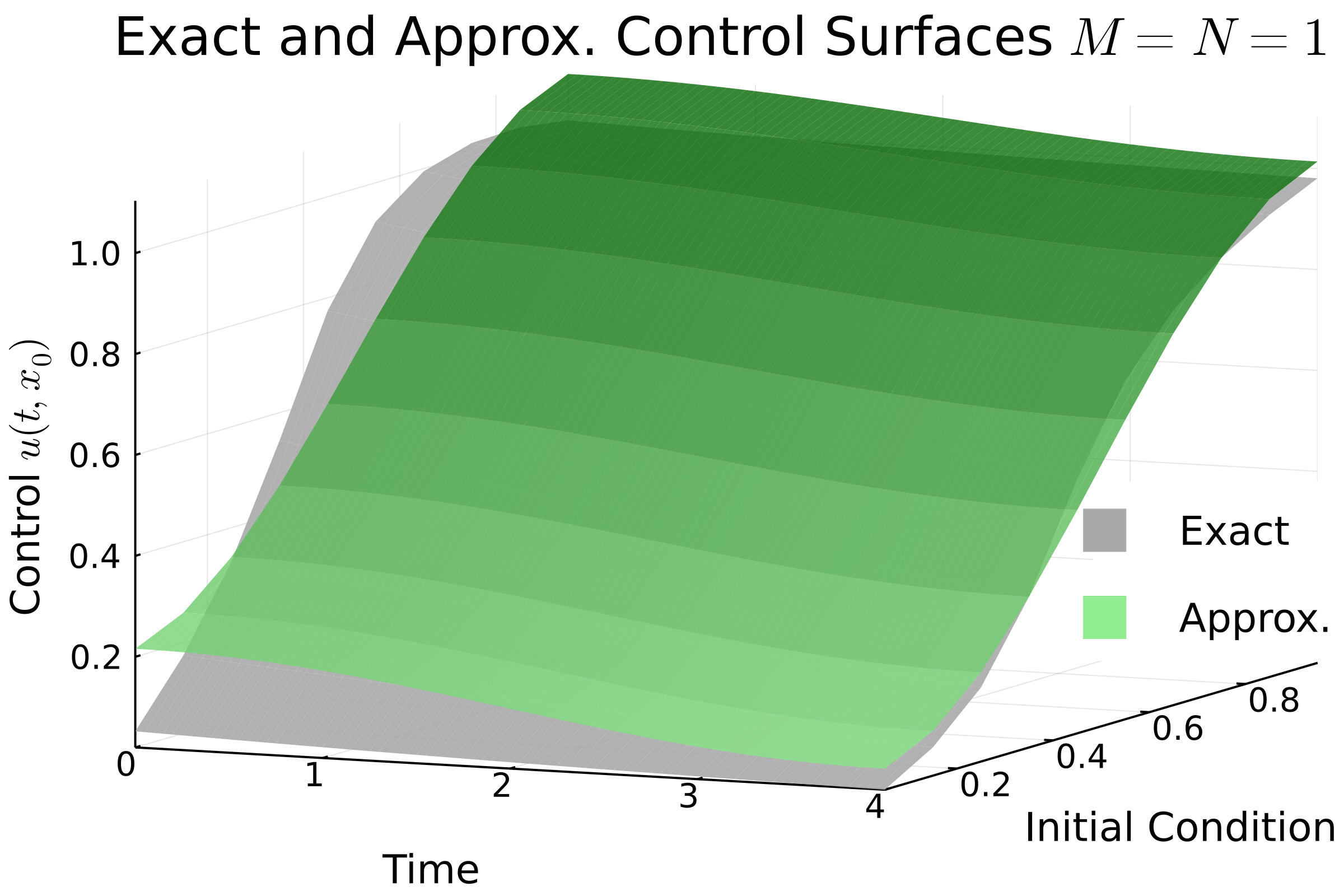} \qquad
\includegraphics[width=0.46\textwidth]{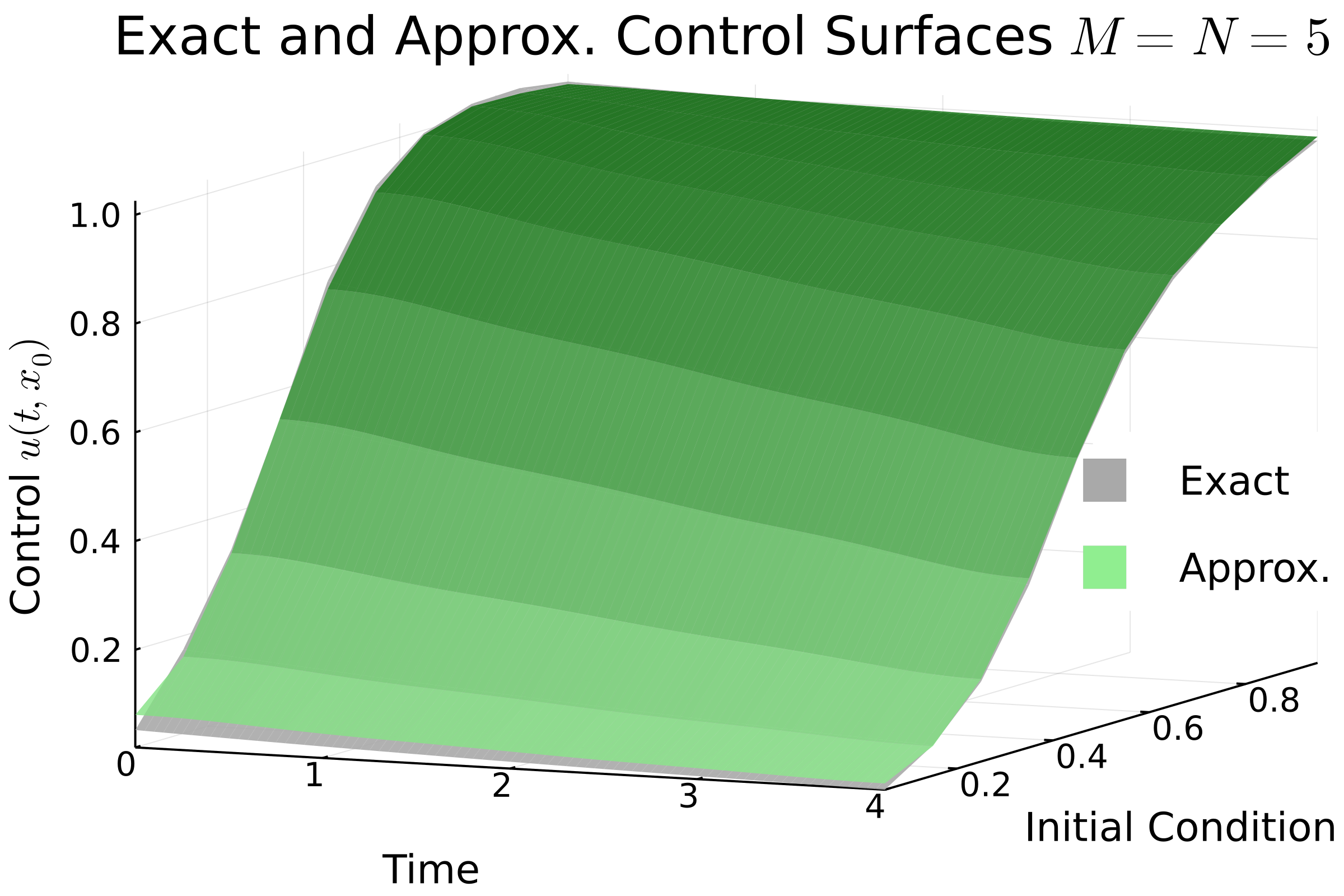}
\caption{(Left) Control surface approximation using a single Fourier term. (Right) Control surface approximation using 5 Fourier terms for each variable.}
\label{figure: controlsurfaces}
\end{figure}

\begin{figure}[ht!]
\centering
\includegraphics[width=0.46\textwidth]{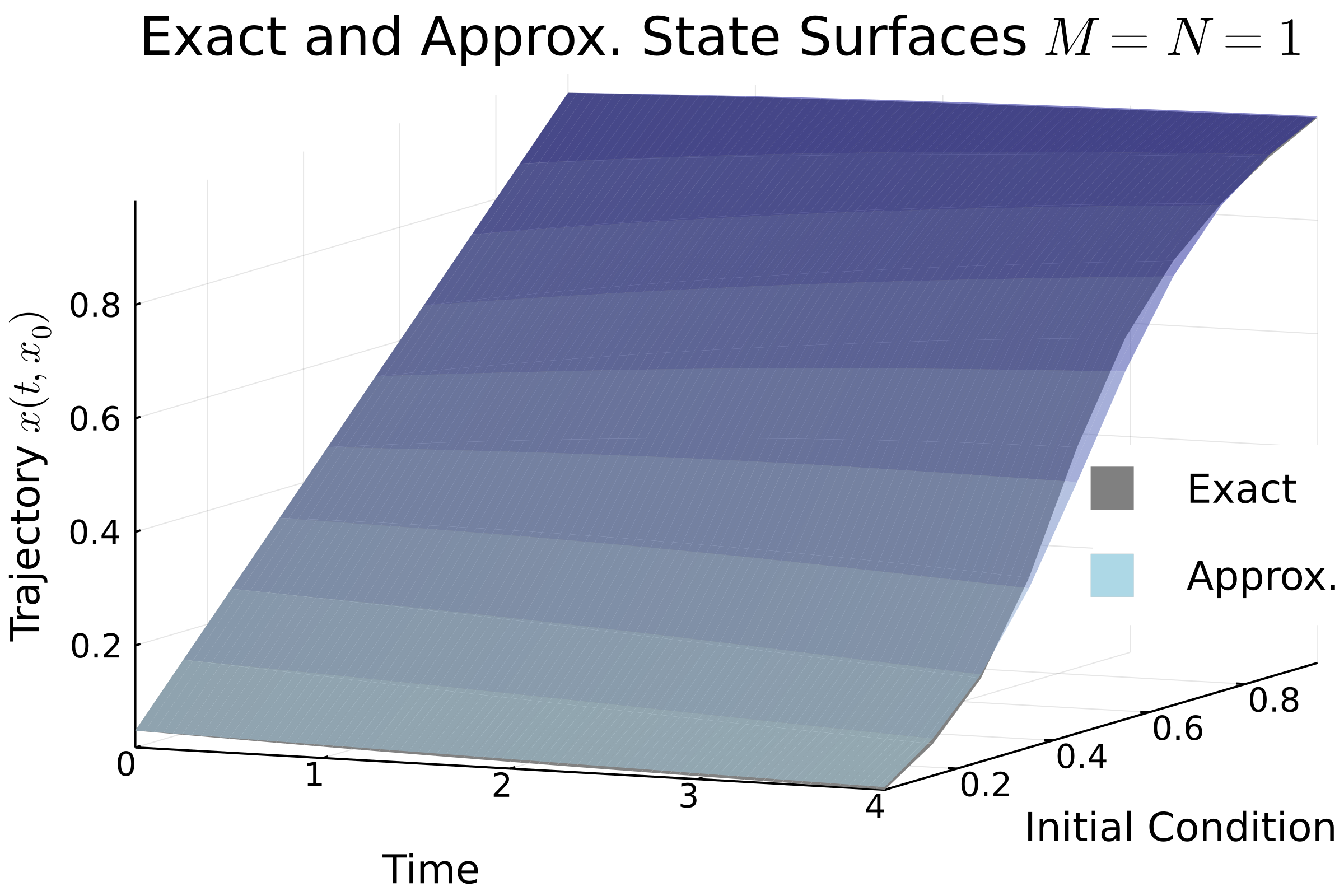} \qquad
\includegraphics[width=0.46\textwidth]{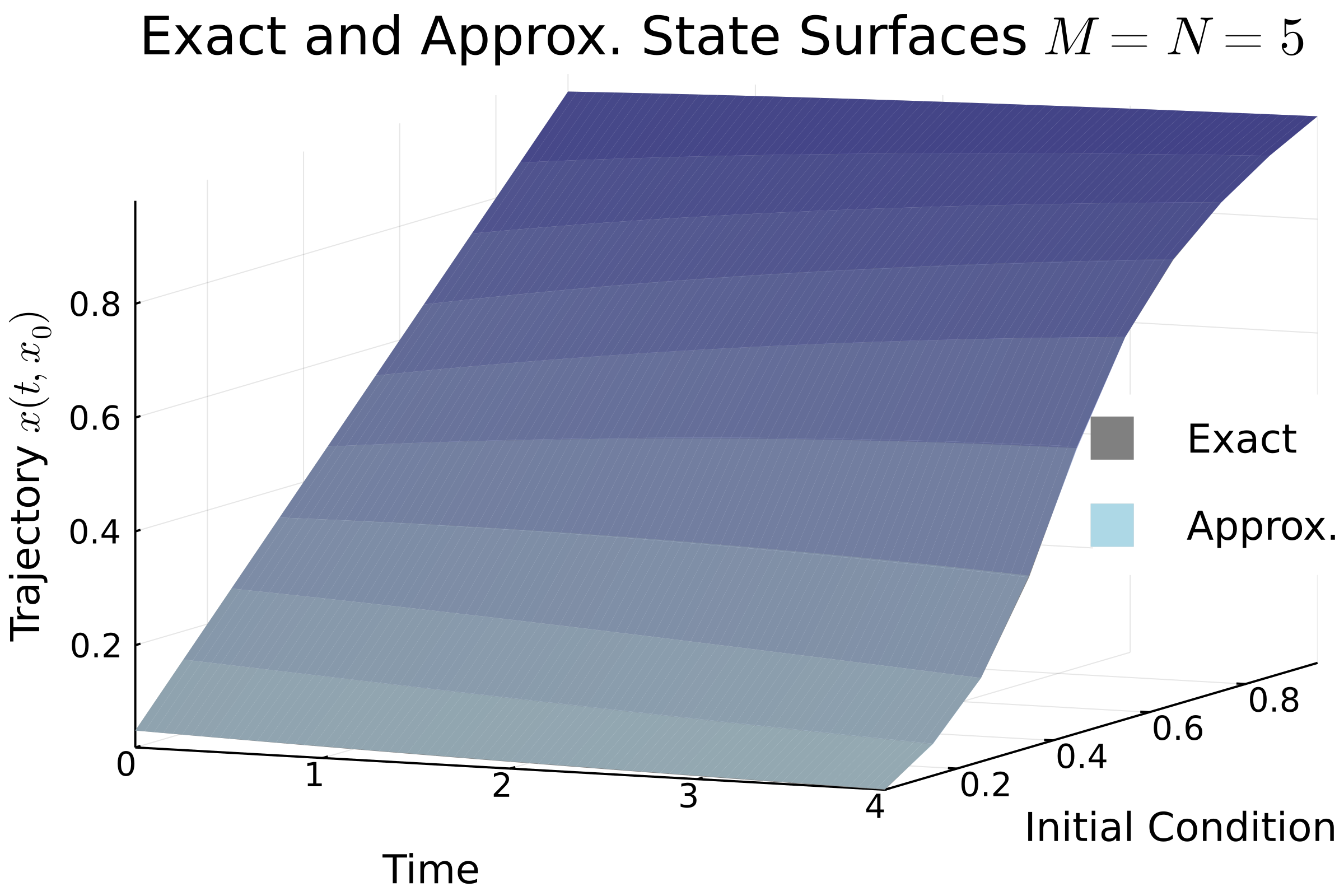}
\caption{(Left) State surface approximation using a single Fourier term. (Right) State surface approximation using 5 Fourier terms for each variable.}
\label{figure: statesurfaces}
\end{figure}

\begin{figure}[ht!]
\centering
\includegraphics[width=0.45\textwidth]{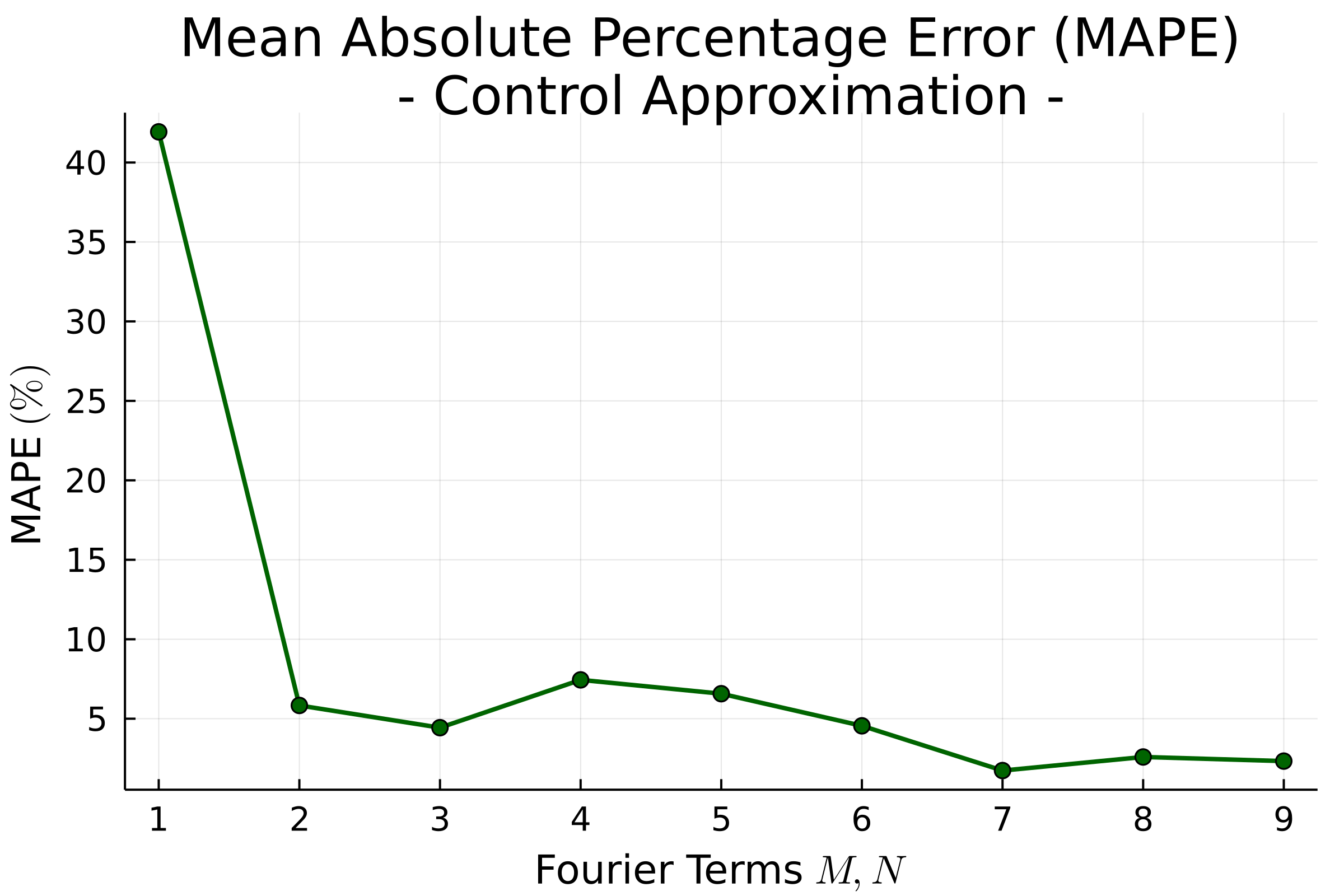} \qquad
\includegraphics[width=0.45\textwidth]{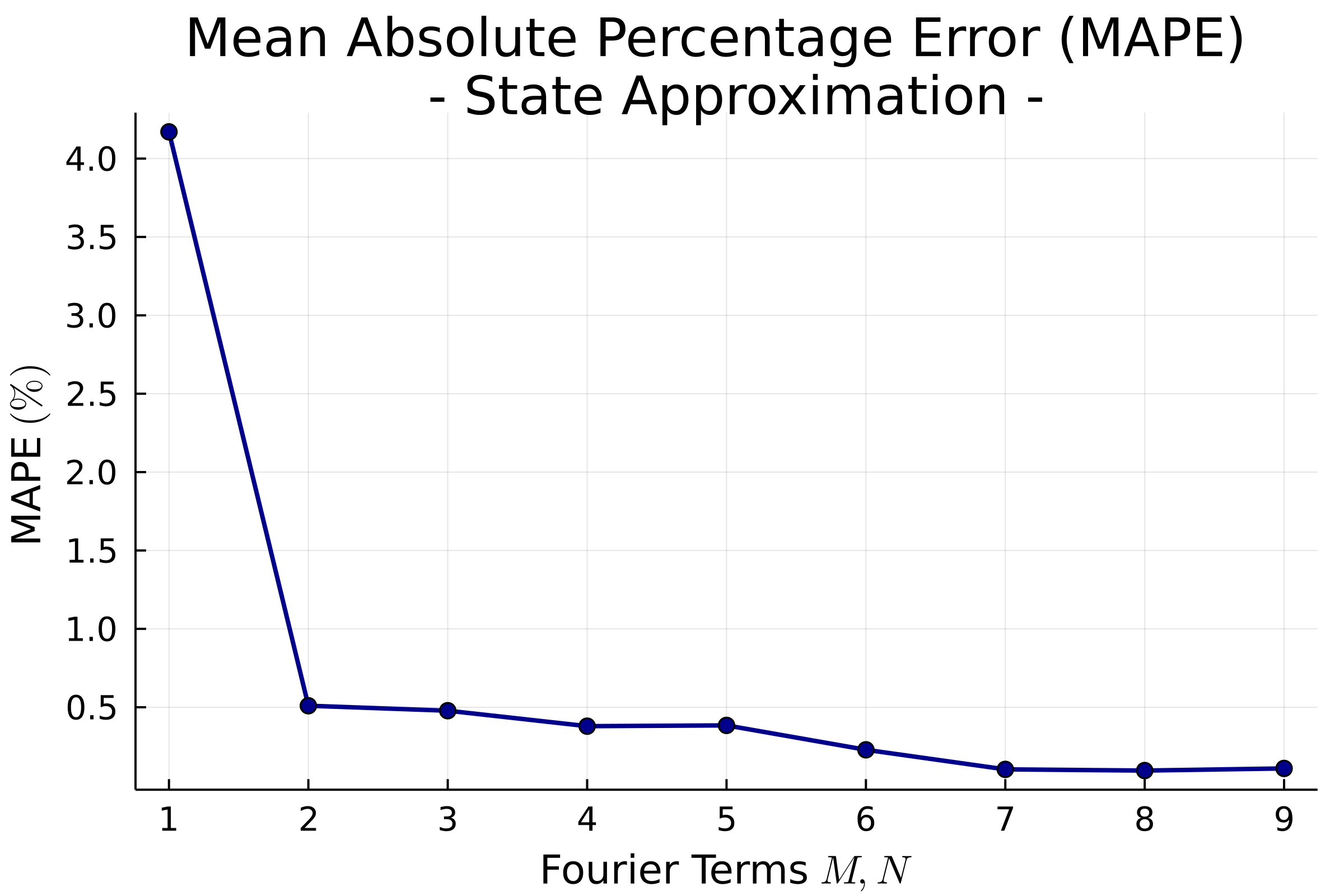}
\caption{Mean Absolute Percentage Error for the control (left) and state (right) approximation for $t \in \{0,0.05,\dots,4.0\}$ and $x_0 \in \{0.05,0.10,\dots, 0.95 \}$, with gradient tolerance $\epsilon = 10^{-4}$.}
\label{figure: MAPEcontrol}
\end{figure}

Results are shown in \cref{figure: controlsurfaces,figure: statesurfaces,figure: MAPEcontrol}. In \cref{figure: controlsurfaces} and \cref{figure: statesurfaces}, two instances of the approximated optimal control and corresponding approximated optimal state surfaces are presented, each considering coefficients $M=N=1$ and $M=N=5$, respectively.     

\begin{figure}[ht!]
\centering
\includegraphics[width=0.65\textwidth]{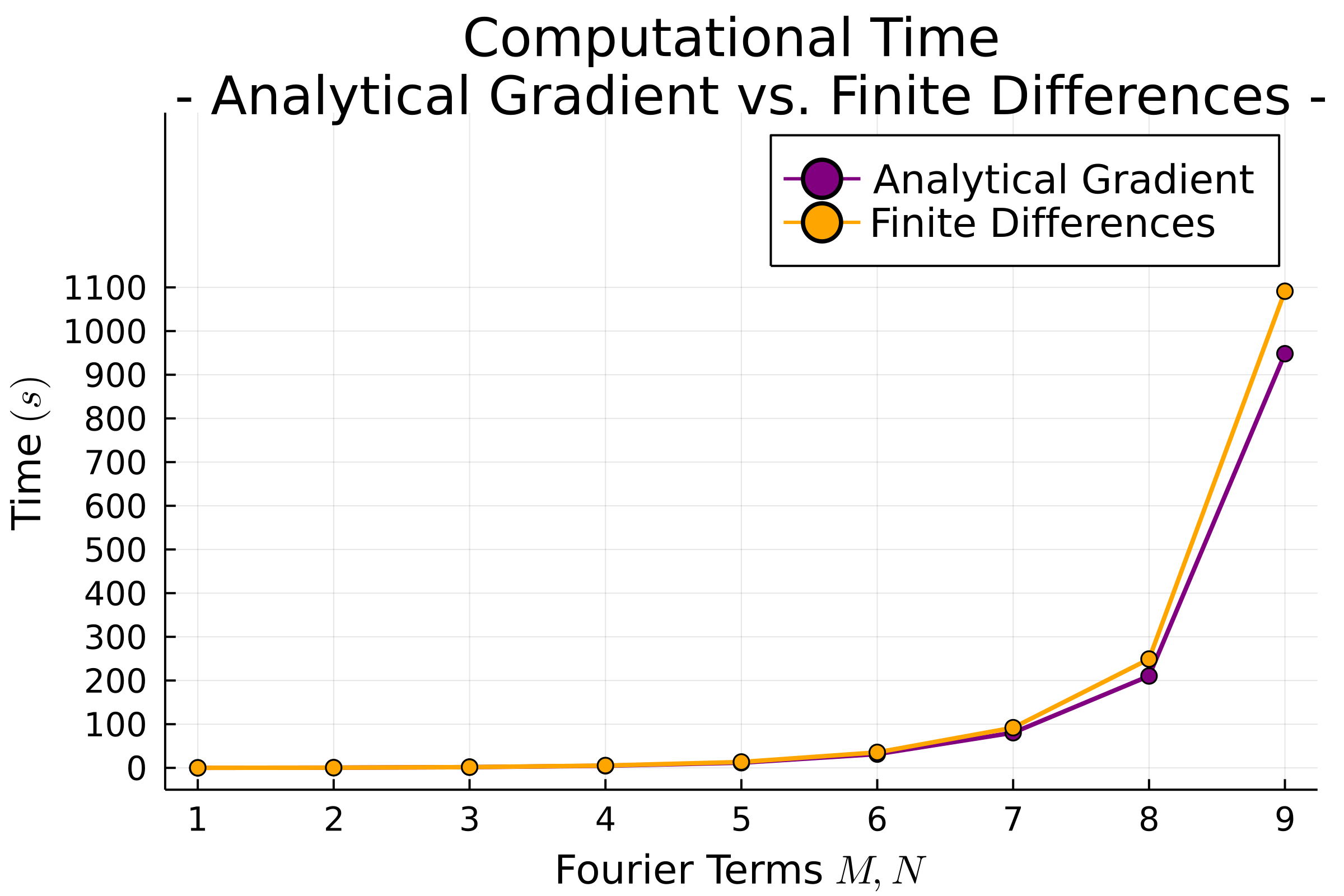} \qquad
\caption{Comparison of computational time of Gradient Descent in \cref{alg:cap}. Analytical gradient (purple line) versus approximated (finite differences) gradient (yellow line).}
\label{figure: GradvsFin}
\end{figure}
\textcolor{black}{\cref{figure: controlsurfaces} shows two approximated surfaces for $M=N=1$ and $M=N=5$, respectively. The control surface on the right is a close approximation to the actual control surface derived from the analytical solution of \cref{eqn:SoftwareMaintenanceProblem}, presenting a mean absolute percentage error (MAPE) of $6.57\%$. State trajectories derived from the control approximations in \cref{figure: controlsurfaces} are shown in \cref{figure: statesurfaces}.}

\textcolor{black}{The state trajectory (\cref{figure: statesurfaces}) shows the non-equilibrium transition behavior caused by control inputs. For the two-strategy game under consideration the two equilibria are pure strategies. When the system is initialized at a mixed strategy the control can be used to switch among these two strategy equilibria or to drive the system toward a specific strategy equilibrium based on the objective function. In this case, the controller is being used to push the system toward an equilibrium that would not naturally occur without a control input. The system does not settle into an equilibrium in the finite time horizon $[0,T]$.
}

\textcolor{black}{
Additionally, \cref{figure: MAPEcontrol} presents MAPEs obtained for other choices for the number of Fourier terms, both for the control (left) and state (right) approximations. In \cite{FG17}, it is shown that the optimal controller is always a decreasing function of time}. This is also clearly shown in \cref{figure: controlsurfaces}, showing that an interpretable theoretical result can be derived from the approximation. \textcolor{black}{Lastly, \cref{figure: GradvsFin} shows the computational time, in seconds, of \cref{alg:cap} when the gradient is provided analytically (purple line) and when the descent procedure uses finite differences approximations for the gradient (yellow line) when solving the problem defined by \cref{eqn:SoftwareMaintenanceProblem}. As expected, as the number of Fourier terms increases, \cref{alg:cap} is seen to benefit from the analytical gradient, justifying the back-propagation operation discussed in \cref{sec:BackPropagation}}.  

\section{Generalizations} \label{sec:Generalization}
The fact that the dynamics in \cref{eqn:xdot} substantially simplified the construction of the back-propagation operator also suggests that this method is generalizable. A fundamental assumption of this paper is that the state evolution is represented by a closed form integrable differential equation. Therefore assume for a general problem that the state constraint in \cref{eqn:OptCon} can be solved as
\begin{equation}
    \textcolor{black}{\mathcal{V}}(x(t)) = \int_0^t w(u(\tau)) \; d\tau + \textcolor{black}{\mathcal{V}}(x_0).
    \label{eqn: separable3}
\end{equation}
Such an example occurs when
\begin{equation}
    \dot{x} = \frac{w(u(t))}{v(x)}
    \label{eqn: separable}
\end{equation}
so that we can rewrite it as
\begin{equation}
    v(x) \; dx = w(u(t)) \; dt,
    \label{eqn: separable2}
\end{equation}
as we have in \cref{eqn:xdot}. However, more complex examples leading to \cref{eqn: separable3} are possible. If we can write
\begin{equation}
    x(t) = \textcolor{black}{\mathcal{V}}^{-1}\left( \int_0^t w\left(u(\tau)\right)d\tau + \textcolor{black}{\mathcal{V}}(x_0) \right)
    \label{eqn: explicit_x}
\end{equation}
where $u(\tau)$, in the context of this paper, is actually $u(\tau; \mathbf{a})$, then as a consequence, this equation for $x(t)$ explicitly is also a function of the Fourier coefficients $a_{mn}$, so one can compute its partial derivative $\frac{\partial x(t; \; \mathbf{a})}{\partial a_{mn}}$. As we saw in \cref{sec:BackPropagation}, the gradient of the objective functional $\nabla J_{a_{mn}}$ will be also computed explicitly. Following the notation in \cref{eqn:OptCon}, we write
\begin{equation}
    \frac{\partial J}{\partial a_{mn}} = \int_0^T \left[ \frac{\partial f}{\partial x}\frac{\partial x}{\partial a_{mn}} + \frac{\partial f}{\partial u}\frac{\partial u}{\partial a_{mn}} \right] \; dt  
\label{eqn: partialJpartiala_mn}
\end{equation}
In \cref{eqn: partialJpartiala_mn}, $\frac{\partial u}{\partial a_{mn}}$ is computed directly from the approximated controller by the proposed truncated Fourier series in \cref{eqn:ufourier}. Similarly, a closed form of $\frac{\partial x}{\partial a_{mn}}$ can also be computed from the obtained $x(t; \mathbf{a})$ in \cref{eqn: explicit_x}. In this manner, $\frac{\partial x}{\partial a_{mn}}$ is
\begin{equation}
    \frac{\partial x}{\partial a_{mn}} = \frac{\partial}{\partial a_{mn}}\left[ \textcolor{black}{\mathcal{V}}^{-1}\left( \int_0^t w\left(u(\tau ; \mathbf{a})\right)d\tau + \textcolor{black}{\mathcal{V}}(x_0) \right) \right].
\end{equation}
This can be used explicitly in the construction of back-propagation operators for more general problems.

\section{Conclusions and Future Directions}
\label{sec:Conclusion}
In this paper we used a Fourier approximation approach to derive a back-propagation operator to numerically construct approximate optimal control surfaces for non-linear control problems arising from two-strategy evolutionary games. Our formulation is phrased as a non-standard feedforward neural network approximation method, justifying the use of the back-propagation construction. We showed empirically that this method works well for constructing high quality control surfaces in an example control problem on evolutionary games dynamics. We also showed that this method is generalizable when the resulting controlled state dynamics are integrable.

In future work we will generalize the evolutionary game problems we consider to more strategies in an attempt to analyze the control problems found in \cite{GF22}. In particular, this will require dealing with non-integrable dynamical systems. The work in \cite{GF22} shows that the general non-linear controller in cyclic games with an odd number of strategies exhibit oscillations whose properties may be elucidated by this Fourier approximation method. In addition, we will consider control problems for evolutionary games on graphs where chaos can emerge \cite{GSB22}. This will provide an interesting case study for the numerical stability of this approach and may shed further light on these newly emerging problems in non-linear dynamics. 

\section*{Acknowledgements}
Portions of G. N.'s work were sponsored by the National Science Foundation under grant DMS-1814876 and CMMI-1932991. Portions of C.G.'s work were supported by the National Science Foundation under grant CMMI-1932991.

\bibliographystyle{apsrev4-1}
\bibliography{cas-refs}

\end{document}